\batchmode
\makeatletter
\def\input@path{{\string"/home/aeftimia/Documents/projects/Discrete Exterior Calculus/\string"/}}
\makeatother
\documentclass[english]{article}
\usepackage[T1]{fontenc}
\usepackage[latin9]{inputenc}
\usepackage{float}
\usepackage{amsmath}
\usepackage{amssymb}
\usepackage{graphicx}
\usepackage{esint}

\makeatletter

\providecommand{\tabularnewline}{\\}

\numberwithin{equation}{section}
\numberwithin{figure}{section}
\numberwithin{table}{section}

\date{}

\@ifundefined{showcaptionsetup}{}{%
 \PassOptionsToPackage{caption=false}{subfig}}
\usepackage{subfig}
\makeatother

\usepackage{babel}
\begin{document}

\title{Kahler: An Implementation of Discrete Exterior Calculus on Hermitian
Manifolds}

\author{Alex Eftimiades%
\thanks{University of Maryland at Baltimore County, aeftimi1@umbc.edu%
}}
\maketitle
\begin{abstract}
This paper details the techniques and algorithms implemented in Kahler,
a Python library that implements discrete exterior calculus on arbitrary
Hermitian manifolds. Borrowing techniques and ideas first implemented
in PyDEC, Kahler provides a uniquely general framework for computation
using discrete exterior calculus. Manifolds can have arbitrary dimension,
topology, bilinear Hermitian metrics, and embedding dimension. Kahler
comes equipped with tools for generating triangular meshes in arbitrary
dimensions with arbitrary topology. Kahler can also generate discrete
sharp operators and implement de Rham maps. Computationally intensive
tasks are automatically parallelized over the number of cores detected.
The program itself is written in Cython--a superset of the Python
language that is translated to C and compiled for extra speed. Kahler
is applied to several example problems: normal modes of a vibrating
membrane, electromagnetic resonance in a cavity, the quantum harmonic
oscillator, and the Dirac-Kahler equation. Convergence is demonstrated
on random meshes.
\end{abstract}

\section{Introduction}

The ideas and techniques used in Kahler%
\footnote{Code available at https://github.com/aeftimia/kahler%
} are based on those pioneered by Bell and Hirani in PyDEC \cite{Bell2012}.
This section will briefly review the core concepts of discrete exterior
calculus. Please refer to \cite{Hirani2005} for a more detailed overview
of discrete exterior calculus, \cite{Bell2012} for its implementation,
and \cite{Abraham1988} for its continuous counterpart.

\subsection{The Discrete Exterior Derivative}

Each $p$-simplex, $\sigma^{p}$, is composed of $p+1$ vertices

\begin{equation}
\sigma^{p}=\left[v_{0},\dots,v_{p}\right]
\end{equation}

In practice, only the indices of these vertices are stored in the
simplices. The boundary of a $p$-simplex is defined as the formal
sum of its $p-1$ dimensional faces \cite{Hirani2005}

\begin{equation}
\partial\sigma^{p}=\underset{i}{\sum}\left(-1\right)^{i}\left[v_{0},\dots,\hat{v_{i}}\dots,v_{p}\right]
\end{equation}

The discrete exterior derivative is defined as the transpose of the
boundary operator \cite{Hirani2005}

\begin{equation}
d=\partial^{T}
\end{equation}

\subsection{The Dual Mesh and the Hodge Star}

The discrete Hodge star maps $p$-simplices, $\sigma^{p}=\left[v_{0},\dots,v_{p}\right]$
to their circumcentric duals \cite{Hirani2005}

\begin{equation}
\star\sigma^{p}=\underset{\sigma^{p+1}\prec\cdots\prec\sigma^{N}}{\sum}s\left(\sigma^{p},\dots,\sigma^{N}\right)\left[c\left(\sigma^{N}\right),\dots,c\left(\sigma^{p}\right)\right]
\end{equation}

, where $c\left(\sigma\right)$ is the circumcenter of $\sigma$,
$\sigma\prec\sigma^{\prime}$ means $\sigma$ is a face of $\sigma^{\prime}$.
$s$ is $1$ when the orientation of $\left[c\left(\sigma^{N}\right),\dots,c\left(\sigma^{p}\right)\right]$
matches the orientation of $\left[\sigma^{N}\backslash\sigma^{N-1},\dots,\sigma^{p+1}\backslash\sigma^{p},c\left(\sigma^{p}\right)\right]$
and is $-1$ otherwise. The discrete Hodge star is then defined as
the following diagonal matrix \cite{Hirani2005,Bell2012}

\begin{equation}
\star_{ij}=\delta_{ij}\frac{\left|\star\sigma_{i}\right|}{\left|\sigma_{j}\right|}
\end{equation}

where $\left|\sigma_{j}\right|$ is the primal volume of the $j^{\mathrm{th}}$
simplex and $\left|\star\sigma_{i}\right|$ is the volume of the dual
cell of the $j^{\mathrm{th}}$ simplex. The codifferential is then
defined as \cite{Hirani2005}

\begin{equation}
\delta_{p}=\left(-1\right)^{p}\star_{p-1}^{-1}d_{p-1}^{T}\star_{p}
\end{equation}

\section{Software Overview}

\subsection{Mesh Generation}

For the purposes of this paper, mesh generation refers to converting
a topological manifold to a set of $N$-simplices. Mesh generation
is therefore separate from embedding the points that form these simplices
in some subset of $\mathbb{C}^{M}$ with $M\ge N$. Kahler comes with
tools that accomplish each of these tasks separately. This section
is concerned with mesh generation.

\subsubsection{Generating Simplicial Meshes in Arbitrary Dimensions}

Kahler comes with two types of triangular mesh generators--one for
asymmetric meshes and one for symmetric meshes as shown below.

\begin{figure}[H]
\begin{centering}
\subfloat{\centering{}\includegraphics[scale=0.3]{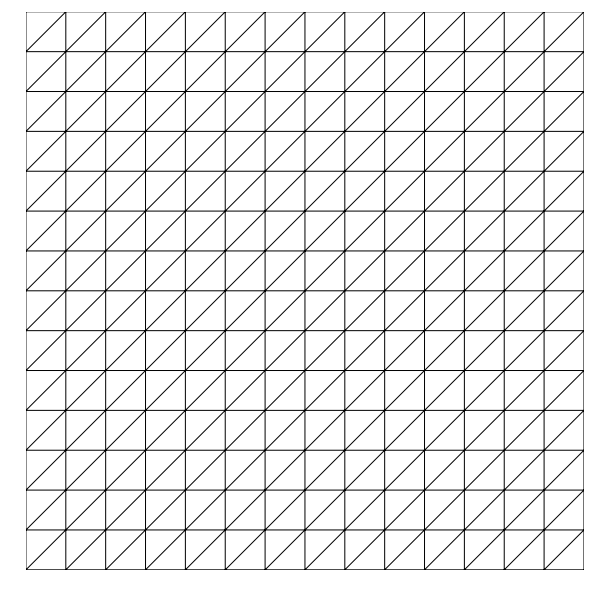}}\subfloat{\centering{}\includegraphics[scale=0.3]{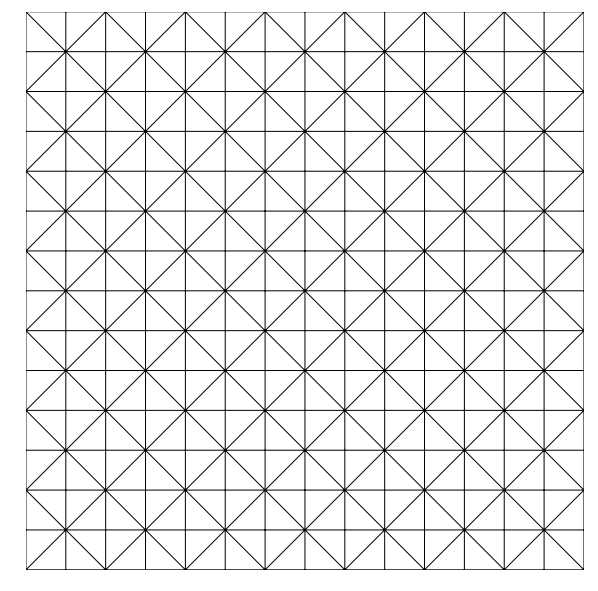}}
\par\end{centering}

\begin{centering}
\subfloat{\centering{}\includegraphics[scale=0.3]{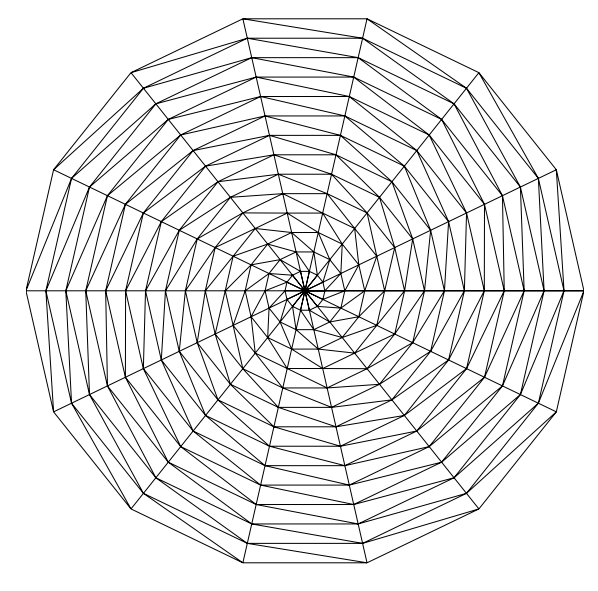}}\subfloat{\centering{}\includegraphics[scale=0.3]{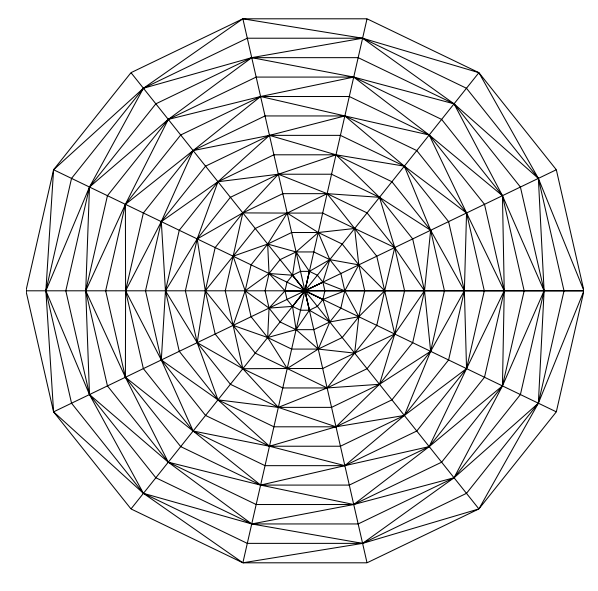}}
\par\end{centering}

\caption{Asymmetric 2D mesh (top left), symmetric 2D mesh (top right), asymmetric
polar mesh (bottom left), symmetric polar mesh (bottom right). }
\end{figure}

Similar meshes can be generated for tori, spheres, etc\@. Meshes
start as an $M_{1}\times M_{2}\times\dots\times M_{N}$ grid of points
created using the grid\_indices() function. Grid indices are stored
as dictionaries that map each index on the grid (grid index) to the
corresponding index of the vertex (vertex index). Simplices are generated
by connecting points to neighboring points. Simplices are stored as
lists of vertex indices in the form of Numpy arrays. Asymmetric grids
connect each point on the grid with a grid index of $\left[n_{1},\dots,n_{i},\dots n_{N}\right]$
to the point with grid index $\left[n_{1},\dots,n_{i}+1,\dots n_{N}\right]$.
Symmetric grids connect every other point on the grid with a grid
index of $\left[n_{1},\dots,n_{i},\dots n_{N}\right]$ to point with
grid indices $\left[n_{1},\dots,n_{i}\pm1,\dots n_{N}\right]$. Let
$\vec{n}_{k}$ be the grid index of the $k^{\mathrm{th}}$ vertex
in the simplex. Simplices are generated by constructing a series grid
indices such that $\vec{n}_{k+1}-\vec{n}_{k}$ is a unit vector orthogonal
to the plane formed by $\left[\vec{n}_{0},\dots,\vec{n}_{k}\right]$.
In practice, only the vectors, $\vec{n}_{k}-\vec{n}_{0}$, are calculated.
These vectors are added to each grid index and the original grid index
dictionary is used to look up the corresponding vertex index--if it
exists.

\subsubsection{Customizing Topology}

Once a list of simplices is created, arbitrary topologies can be created
by joining one or more vertices together. For example, a torus is
formed by joining vertices at opposite ends of the grid in question.
This processes is hereby referred to as stitching. Stitches are stored
as a dictionary that maps vertex indices to the corresponding vertex
indices they are joined to. Kahler comes equipped with a function,
pbc\_stitches(), that will create the stitches necessary for periodic
boundary conditions in one or more directions given a set of grid
indices.

\subsubsection{Singular Hodge Stars}

The type of mesh discussed in this section creates right $N$-simplices.
This means that the circumcenter of any $p$-simplex lies on its $p-1$
dimensional hypotenuse. The leads to dual cells with zero volume,
which means that the inverse of the discrete Hodge star operator cannot
be defined on those cells. This problem was remedied by adding a small
but nonzero real number to each entry of the metric at each point.
It is worth noting that this problem does not occure when calculating
the Laplace-Beltrami operator for $0$-forms because the only dual
volumes that need to be inverted are those of the $0$-simplices.
These are generally nonzero even for right simplices.

\subsection{The de Rham Map}

A de Rham map is a mapping from a continuous differential form to
a discrete one. Given a continuous $p$-form, $\alpha$, a de Rham
map, $R$, is defined as follows

\begin{equation}
R:\alpha\Rightarrow\left(\sigma\rightarrow\int_{\sigma}\alpha\right)
\end{equation}

where $\sigma$ is a $p$-simplex. When $\alpha$ is a continuous
tensor field, the de Rham map is simply a zeroth order polynomial
interpolation of $\alpha$ over the simplicial complex in question.
Kahler uses a generalized trapezoidal rule to compute de Rham maps.
Points are sampled evenly throughout the interior of each $N$-simplex.
The interior of each $\left(N-k\right)$-face is given a weight of
$1/2^{k}$.

\subsection{Computing Circumcenters}

Circumcenters are computed by solving a matrix equation. Given a (locally)
flat Hermitian metric, $h$, a circumcenter can be defined as the
linear combination of vertices, $c^{k}=\alpha^{I}v_{I}^{k}$, such
that $\overline{\left(c^{j}-v_{I}^{j}\right)}h_{jk}\left(c^{k}-v_{I}^{k}\right)=\overline{\left(c^{j}-v_{0}^{j}\right)}h_{jk}\left(c^{k}-v_{0}^{k}\right)$
for each $1\le I\le p$, for a given $p$-simplex. Simplifying yields

\begin{equation}
2\Re\left(\overline{v}_{J}^{j}h_{jk}\left(v_{I}^{k}-v_{0}^{k}\right)\right)\alpha^{J}=\overline{v}_{I}^{j}h_{jk}v_{I}^{k}-\overline{v}_{0}^{j}h_{jk}v_{0}^{k}
\end{equation}

This, combined with the constraint that $\sum\alpha^{I}=1$, uniquely
determines the circumcenter.

\subsection{Orienting Dual Volumes}

Let the circumcenter of a $p$-simplex, $\sigma^{p}$ be $c\left(\sigma^{p}\right)$.
A dual simplex of $\sigma^{p}$ is a simplex formed by $\left[c\left(\sigma^{N}\right),\dots,c\left(\sigma^{p}\right)\right]$
where $\sigma^{p}\prec\cdots\prec\sigma^{N}$. Let $V_{I}=c\left(\sigma^{p+I}\right)-c\left(\sigma^{p}\right)$
for $1\le I\le N$ and $\sigma^{p}\prec\cdots\prec\sigma^{p+I}\prec\cdots\prec\sigma^{N}$.
Furthermore, let $U_{I}=\sigma^{p+I}\backslash\sigma^{p+I-1}-c\left(\sigma^{p}\right)$.
Then the sign associated with a dual simplex, $\left[c\left(\sigma^{N}\right),\dots,c\left(\sigma^{p}\right)\right]$,
is given by

\begin{equation}
\mathrm{sgn}\left(\Re\left(\det\left(\overline{V}^{i}h_{ij}U^{j}\right)\right)\right)
\end{equation}

The metric, $h$, is the metric associated with $\sigma^{N}$.

\subsection{Computing Barycentric Differentials}

Barycentric differentials, $X$, are computed in a similar manner
to the technique Bell and Hirani used for PyDEC. In Euclidean space
\cite{Bell2012}

\begin{equation}
X^{T}=\left(VV^{T}\right)^{-1}V
\end{equation}

with $V_{I}^{i}=v_{I}^{i}-v_{0}^{i}$ for $I\in\left\{ 1,\dots,N\right\} $.
This determines $N$ of the $N+1$ barycentric differentials. The
last is determined from the constraint $\sum X^{I}=0$. In a space
endowed with a Hermitian metric, $h$, this is easily modified accordingly.

\begin{equation}
X^{\dagger}=\left(V^{i}h_{ij}V^{\dagger j}\right)^{-1}V^{k}h_{kl}
\end{equation}

\subsection{The Discrete Sharp Operator}

The discrete sharp operator maps discrete $p$-forms to antisymmetric
tensors located at the barycenters of each $N$-simplex. Kahler wedges
barycentric differentials to form antisymmetric tensors. The wedged
differentials associated with each $p$-face of each $N$-simplex
are averaged to give the tensor associated with the given $N$-simplex.

\begin{equation}
\alpha^{\sharp}\left(\sigma^{N}\right)=\frac{1}{A}\underset{i_{k},j=0}{\overset{N}{\sum}}\alpha\left(\left[\sigma_{i_{1}}^{N},\dots,\sigma_{i_{p}}^{N}\right]\right)\left(-1\right)^{j}d\lambda_{i_{1}}\wedge\dots\wedge\hat{d\lambda_{i_{j}}}\wedge\dots\wedge d\lambda_{i_{p}}
\end{equation}

where $\sigma^{N}$ is an $N$-simplex, $\sigma_{i}^{N}$ is the $i^{\mathrm{th}}$
vertex of $\sigma^{N}$ such that $\left[\sigma_{i_{1}}^{N},\dots,\sigma_{i_{p}}^{N}\right]$
forms a $p$-simplex. $d\lambda_{i}$ is the $i^{\mathrm{th}}$ barycentric
differential at $\sigma^{N}$. $A$ is the total number of terms in
the sum, given by $A={n+1 \choose p+1}{p+1 \choose p}=\frac{\left(n+1\right)!}{p!\left(n+1-p\right)!}$.
In practice, the discrete sharp operator is stored as a compressed
matrix for speed and memory efficiency. When it is dotted with a $p$-form,
it gives a flattened list of antisymmetric tensors. That is, the result
is a one dimensional array that can be reshaped into a list of the
desired antisymmetric tensors.

\section{Examples}

\subsection{Normal Modes of a Vibrating Membrane}

A vibrating membrane can be modeled as a scalar field, $\phi$, that
satisfies the condition $-\delta d\phi=\lambda\phi$. When the membrane
is held taught at the boundary, it obeys the Dirichlet condition $\phi=0$
at the boundary. When the membrane is allowed to vibrate freely at
the boundary (as is the case for a Chladni plate), $\delta\phi=0$
at the boundary. Operating $\delta$ on both sides of the original
eigenvalue equation and making use of $\delta^{2}=0$ yields $-\delta\delta d\phi=0=\delta\lambda\phi$.
So this condition is automatically satisfied by solving the eigenvalue
problem with the usual Laplace-Beltrami operator.

In each case, a $20\times20$ grid was used over the interval $\left[0,\pi\right]\times\left[0,\pi\right]$.
For visualization, the discrete $0$-forms were sharpened and interpolated
over a fine grid.

\begin{figure}[H]
\begin{centering}
\subfloat{\centering{}\includegraphics[scale=0.3]{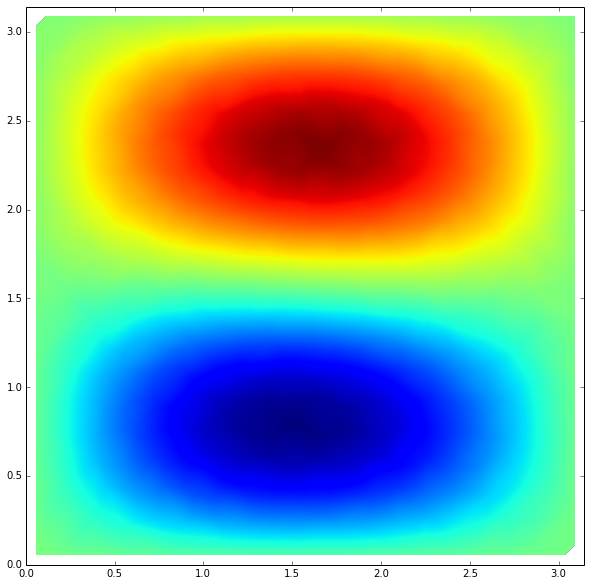}}\subfloat{\centering{}\includegraphics[scale=0.3]{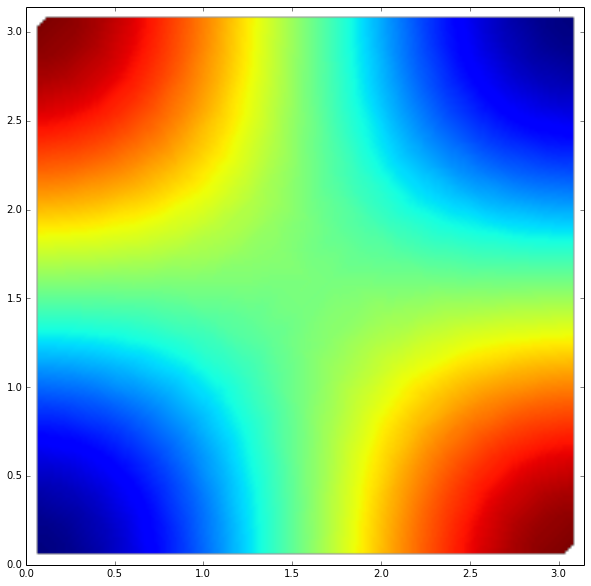}}
\par\end{centering}

\caption{The third normal mode of a drum with Dirichlet boundary conditions
(left) and Neumann boundary conditions (right).}

\end{figure}

\begin{table}[H]
\begin{centering}
\begin{tabular}{|c|c|c|}
\hline 
Analytic &
Numerical &
\%Error\tabularnewline
\hline 
\hline 
2 &
1.99 &
0.426\tabularnewline
\hline 
5 &
4.95 &
1.03\tabularnewline
\hline 
5 &
4.95 &
0.914\tabularnewline
\hline 
8 &
7.91 &
1.11\tabularnewline
\hline 
10 &
9.80 &
2.05\tabularnewline
\hline 
10 &
9.80 &
2.05\tabularnewline
\hline 
13 &
12.8 &
1.95\tabularnewline
\hline 
13 &
12.8 &
1.82\tabularnewline
\hline 
17 &
16.4 &
3.59\tabularnewline
\hline 
17 &
16.4 &
3.58\tabularnewline
\hline 
\end{tabular}
\par\end{centering}

\caption{The First 10 eigenfrequencies of a drum with Dirichlet boundary conditions.}

\end{table}

\begin{table}[H]
\begin{centering}
\begin{tabular}{|c|c|c|}
\hline 
Analytic &
Numerical &
\%Error\tabularnewline
\hline 
\hline 
1 &
0.99 &
0.587\tabularnewline
\hline 
1 &
1.00 &
0.266\tabularnewline
\hline 
2 &
1.99 &
0.427\tabularnewline
\hline 
4 &
3.96 &
1.11\tabularnewline
\hline 
4 &
3.96 &
1.11\tabularnewline
\hline 
5 &
4.95 &
1.09\tabularnewline
\hline 
5 &
4.96 &
0.851\tabularnewline
\hline 
8 &
7.91 &
1.11\tabularnewline
\hline 
9 &
8.80 &
2.25\tabularnewline
\hline 
9 &
8.80 &
2.21\tabularnewline
\hline 
\end{tabular}
\par\end{centering}

\caption{The first 10 eigenfrequencies mode of a square drum with Neumann boundary
conditions.}
\end{table}

\subsection{Electromagnetic Resonance in a Cavity}

Electromagnetic waves in a cavity can be written in terms of the electric
field, $\vec{E}$, or the magnetizing field, $\vec{H}$ \cite{jackson1975classical}.
The governing equations are the same for both, but the boundary conditions
are different. The electric field obeys the Dirichlet condition that
at the boundary, the component of $\vec{E}$ parallel to the boundary
is zero. Let $E=\vec{E}^{\flat}$ and $H=\vec{H}^{\flat}$. When read
in terms of differential forms, $E=0$ at the boundary. At the boundary,
the $\vec{H}$ field is zero perpendicular to the boundary. When read
in terms of differential forms, $\delta H=0$ at the boundary. Since
$\delta dH=\lambda H$, operating with the codifferential on both
sides and making use of $\delta^{2}=0$ yields $\delta\delta dH=0=\delta\lambda H$.
So this condition is automatically satisfied by solving the eigenvalue
problem with the usual Laplace-Beltrami operator.

In each case, a $20\times20$ grid was used over the interval $\left[0,\pi\right]\times\left[0,\pi\right]$.
For visualization, the discrete $1$-forms were sharpened and interpolated
over a fine grid.

\begin{figure}[H]
\begin{centering}
\subfloat{\centering{}\includegraphics[scale=0.3]{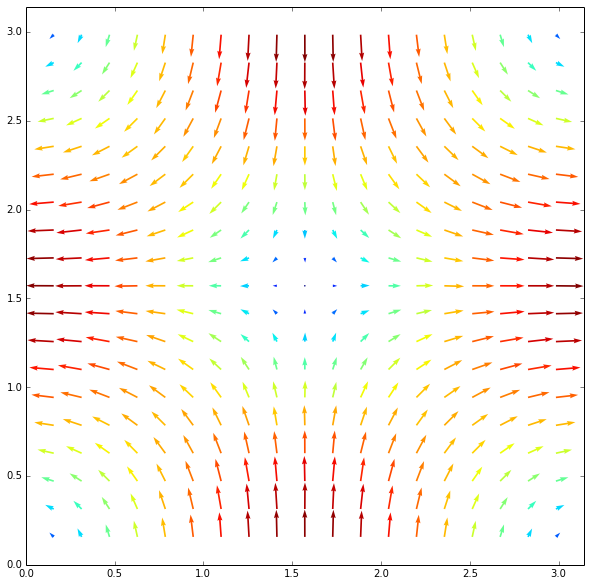}}\subfloat{\centering{}\includegraphics[scale=0.3]{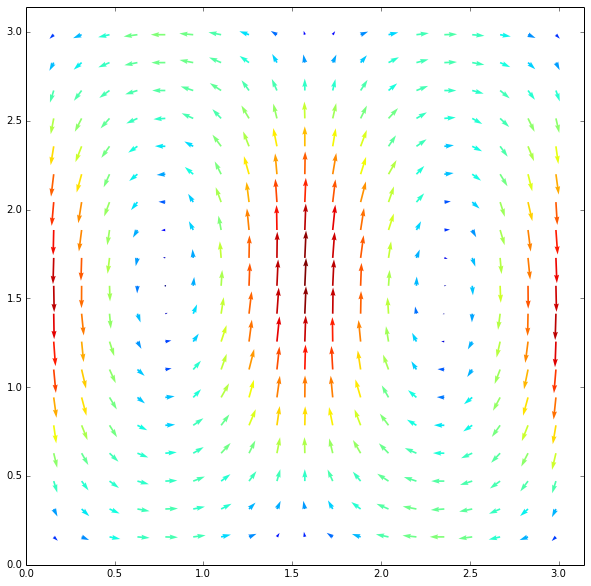}}
\par\end{centering}

\caption{The third normal mode of electric field (left) and magnetic field
(right) in a resonant cavity.}
\end{figure}

\begin{table}[H]
\begin{centering}
\begin{tabular}{|c|c|c|}
\hline 
Analytic &
Numerical &
\%Error\tabularnewline
\hline 
\hline 
1 &
0.996 &
0.386\tabularnewline
\hline 
1 &
1.00 &
0.080\tabularnewline
\hline 
2 &
2.00 &
0.233\tabularnewline
\hline 
4 &
3.96 &
1.06\tabularnewline
\hline 
4 &
3.97 &
0.763\tabularnewline
\hline 
5 &
4.96 &
0.869\tabularnewline
\hline 
5 &
4.97 &
0.685\tabularnewline
\hline 
8 &
7.93 &
0.913\tabularnewline
\hline 
9 &
8.80 &
2.19\tabularnewline
\hline 
9 &
8.83 &
1.89\tabularnewline
\hline 
\end{tabular}
\par\end{centering}

\caption{The first 10 eigenfrequencies of an electric field in a resonant cavity.}
\end{table}

\begin{table}[H]
\begin{centering}
\begin{tabular}{|c|c|c|}
\hline 
Analytic &
Numerical &
\%Error\tabularnewline
\hline 
\hline 
2 &
2.00 &
0.233\tabularnewline
\hline 
5 &
4.96 &
0.868\tabularnewline
\hline 
5 &
4.97 &
0.687\tabularnewline
\hline 
8 &
7.93 &
0.913\tabularnewline
\hline 
10 &
9.80 &
1.98\tabularnewline
\hline 
10 &
9.83 &
1.74\tabularnewline
\hline 
13 &
12.8 &
1.75\tabularnewline
\hline 
13 &
12.8 &
1.64\tabularnewline
\hline 
17 &
16.4 &
3.53\tabularnewline
\hline 
17 &
16.4 &
3.27\tabularnewline
\hline 
\end{tabular}
\par\end{centering}

\caption{The first 10 eigenfrequencies of a of magnetic field in a resonant
cavity.}
\end{table}

\subsection{The Quantum Harmonic Oscillator}

An $M\times\dots\times M$ asymmetric Cartesian grid was used to solve
the quantum harmonic oscillator problem \cite{Shankar1994}.

\begin{equation}
\frac{1}{2}\nabla^{2}\psi+\frac{1}{2}x^{2}\psi=E\psi
\end{equation}

A de Rham map was used to map the potential to the barycenters of
the $N$-simplices. Call this representation of the potential $V^{\sharp}$.
A least squares fit was used to derive the $0$-form representation
of the potential, $V$ (effectively inverting the sharp operator for
$0$-forms).

\paragraph{One Dimension}

A $20$ point grid on the interval $\left[-\pi,\pi\right]$ was used
for the one dimensional harmonic oscillator.

\begin{table}[H]
\begin{centering}
\begin{tabular}{|c|c|c|}
\hline 
Analytic &
Numerical &
\%Error\tabularnewline
\hline 
\hline 
0.5 &
0.483 &
3.37\tabularnewline
\hline 
1.5 &
1.47 &
1.80\tabularnewline
\hline 
2.5 &
2.47 &
1.35\tabularnewline
\hline 
3.5 &
3.49 &
0.164\tabularnewline
\hline 
4.5 &
4.60 &
2.19\tabularnewline
\hline 
5.5 &
5.81 &
5.64\tabularnewline
\hline 
6.5 &
7.13 &
9.73\tabularnewline
\hline 
7.5 &
8.55 &
14.1\tabularnewline
\hline 
8.5 &
10.0 &
17.8\tabularnewline
\hline 
9.5 &
11.5 &
21.3\tabularnewline
\hline 
\end{tabular}
\par\end{centering}

\caption{The first 10 eigenfrequencies of a one dimensional quantum harmonic
oscillator.}
\end{table}

\paragraph{Two Dimensions}

A $20\times20$ point grid on the interval $\left[-\pi,\pi\right]\times\left[-\pi,\pi\right]$
was used for the two dimensional harmonic oscillator.

\begin{table}[H]
\begin{centering}
\begin{tabular}{|c|c|c|}
\hline 
Analytic &
Numerical &
\%Error\tabularnewline
\hline 
\hline 
1 &
0.97 &
2.76\tabularnewline
\hline 
2 &
1.96 &
1.89\tabularnewline
\hline 
2 &
1.96 &
1.89\tabularnewline
\hline 
3 &
2.95 &
1.60\tabularnewline
\hline 
3 &
2.96 &
1.48\tabularnewline
\hline 
3 &
2.96 &
1.481\tabularnewline
\hline 
4 &
3.95 &
1.37\tabularnewline
\hline 
4 &
3.95 &
1.365\tabularnewline
\hline 
4 &
3.98 &
0.407\tabularnewline
\hline 
4 &
3.98 &
0.407\tabularnewline
\hline 
\end{tabular}
\par\end{centering}

\caption{The first 10 eigenfrequencies of a two dimensional quantum harmonic
oscillator.}
\end{table}

\paragraph{Three Dimensions}

A $20\times20\times20$ point grid on the interval $\left[-\pi,\pi\right]\times\left[-\pi,\pi\right]\times\left[-\pi,\pi\right]$
was used for the three dimensional harmonic oscillator.

\begin{table}[H]
\begin{centering}
\begin{tabular}{|c|c|c|}
\hline 
Analytic &
Numerical &
\%Error\tabularnewline
\hline 
\hline 
1.5 &
1.46 &
2.45\tabularnewline
\hline 
2.5 &
2.45 &
1.88\tabularnewline
\hline 
2.5 &
2.45 &
1.88\tabularnewline
\hline 
2.5 &
2.45 &
1.88\tabularnewline
\hline 
3.5 &
3.44 &
1.63\tabularnewline
\hline 
3.5 &
3.44 &
1.63\tabularnewline
\hline 
3.5 &
3.44 &
1.63\tabularnewline
\hline 
3.5 &
3.45 &
1.53\tabularnewline
\hline 
3.5 &
3.45 &
1.53\tabularnewline
\hline 
3.5 &
3.45 &
1.53\tabularnewline
\hline 
\end{tabular}
\par\end{centering}

\caption{The first 10 eigenfrequencies of a three dimensional quantum harmonic
oscillator.}
\end{table}

\subsection{The Dirac-Kahler Equation}

The Dirac-Kahler operator, $d+\delta$, operates on a formal sum of
differential forms of different dimension \cite{Red'kov2011}. A $20\times20$
point grid on the interval $\left[0,\pi\right]\times\left[0,\pi\right]$
was used for the two dimensional Dirac-Kahler equation with a Euclidean
metric. Dirichlet boundary conditions were applied to each of the
three $p$-forms.

\begin{table}[H]
\begin{centering}
\begin{tabular}{|c|c|c|}
\hline 
Analytic &
Numerical &
\%Error\tabularnewline
\hline 
\hline 
1 &
0.997 &
0.272\tabularnewline
\hline 
1 &
0.999 &
0.104\tabularnewline
\hline 
1.41 &
1.41 &
0.189\tabularnewline
\hline 
2 &
1.99 &
0.554\tabularnewline
\hline 
2 &
1.99 &
0.503\tabularnewline
\hline 
2.24 &
2.22 &
0.523\tabularnewline
\hline 
2.24 &
2.23 &
0.400\tabularnewline
\hline 
2.28 &
2.81 &
0.530\tabularnewline
\hline 
3 &
2.97 &
1.13\tabularnewline
\hline 
3 &
2.97 &
1.07\tabularnewline
\hline 
\end{tabular}
\par\end{centering}

\caption{The first 10 eigenvalues of the two dimensional Dirac-Kahler operator
with a Euclidean metric.}
\end{table}

\subsection{Random Meshes}

Some care must be taken when generating random meshes. Points cannot
be too close or the resulting simplices will have zero volume. Furthermore,
if an $N$ dimensional rectangular domain is to be covered randomly,
all $p$ dimensional planes must also contain at least one point.
The first requirement was satisfied by resampling points until all
points were no closer than a cutoff distance, $R$. The second requirement
can be satisfied by breaking up the point sampling into sampling at
all corners, then all lines, two dimensional faces, etc. All the while,
points are resampled until none are closer than $R$ to each other
and all previously sampled points. In practice, $R=\frac{1}{2}\frac{1}{M^{1/N}}$,
where $M$ is the number of points on the mesh. Finally, the resulting
points are Delaunay triangulated. Kahler comes with a function, random\_mesh(),
that can create arbitrary dimensional random meshes with periodicity
in an any number and combination of directions. For testing purposes,
a two dimensional mesh was generated for different numbers of points

\begin{figure}[H]
\begin{centering}
\subfloat{\centering{}\includegraphics[scale=0.3]{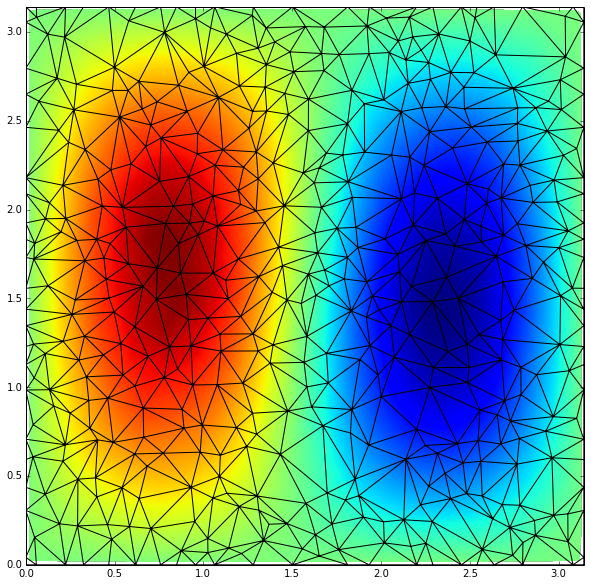}}\subfloat{\centering{}\includegraphics[scale=0.3]{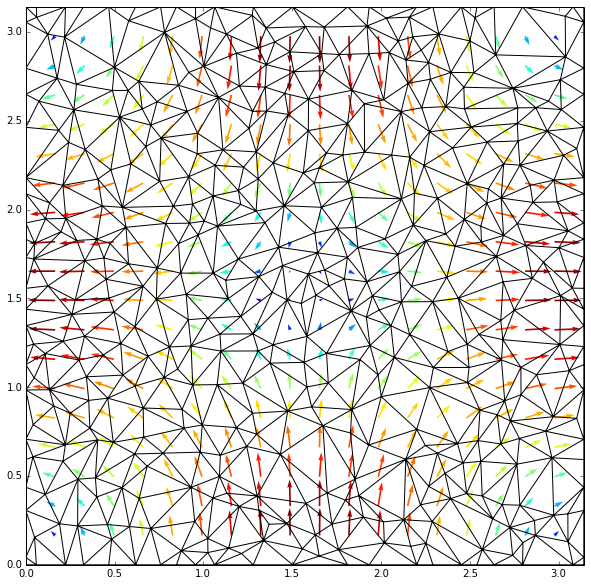}}
\par\end{centering}

\caption{The third normal mode of vibrating drum with Dirichlet boundary conditions
(left) and electric field (right) in a resonant cavity.}
\end{figure}

\begin{table}[H]
\begin{centering}
\begin{tabular}{|c|c|c|}
\hline 
Analytic &
Numerical &
\%Error\tabularnewline
\hline 
\hline 
2 &
1.99 &
0.391\tabularnewline
\hline 
5 &
4.95 &
1.04\tabularnewline
\hline 
5 &
4.97 &
0.598\tabularnewline
\hline 
8 &
7.86 &
1.76\tabularnewline
\hline 
10 &
9.78 &
2.16\tabularnewline
\hline 
10 &
9.85 &
1.48\tabularnewline
\hline 
13 &
12.6 &
3.29\tabularnewline
\hline 
13 &
12.7 &
2.44\tabularnewline
\hline 
17 &
16.3 &
3.83\tabularnewline
\hline 
17 &
16.4 &
3.41\tabularnewline
\hline 
\end{tabular}
\par\end{centering}

\caption{The first 10 eigenfrequencies of a vibrating drum with Dirichlet boundary
conditions.}
\end{table}

\begin{table}[H]
\begin{centering}
\begin{tabular}{|c|c|c|}
\hline 
Analytic &
Numerical &
\%Error\tabularnewline
\hline 
\hline 
1 &
0.99 &
0.245\tabularnewline
\hline 
1 &
1.00 &
0.178\tabularnewline
\hline 
2 &
1.99 &
0.466\tabularnewline
\hline 
4 &
3.96 &
0.778\tabularnewline
\hline 
4 &
3.96 &
0.518\tabularnewline
\hline 
5 &
4.95 &
1.37\tabularnewline
\hline 
5 &
4.96 &
0.858\tabularnewline
\hline 
8 &
7.91 &
1.19\tabularnewline
\hline 
9 &
8.80 &
1.48\tabularnewline
\hline 
9 &
8.80 &
0.989\tabularnewline
\hline 
\end{tabular}
\par\end{centering}

\caption{The first 10 eigenfrequencies of an electric field in a resonant cavity.}
\end{table}

\subsection{Convergence}

Convergence was analyzed for square membrane with Dirichlet boundary
conditions on a random mesh. The average percent error of the first
10 eigenfrequencies of this membrane was computed for increasing numbers
of vertices. The results were repeated five times and averaged.

\begin{figure}[H]
\begin{centering}
\includegraphics[scale=0.3]{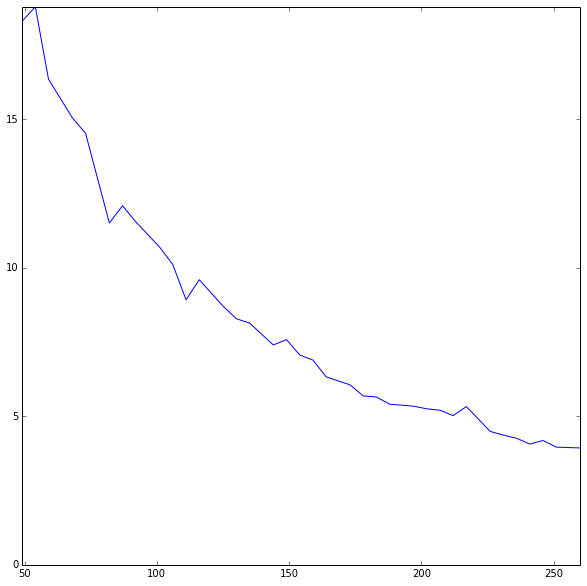}
\par\end{centering}

\caption{The mean of percent errors of the first 10 eigenfrequencies of a square
membrane with Dirichlet boundary conditions vs the number of points
on a random mesh.}
\end{figure}

\section*{Acknowledgments}

I thank Dr. Muruhan Rathinam for helping me learn the foundations
of discrete exterior calculus.

\bibliographystyle{plain}
\bibliography{11_home_aeftimia_Documents_projects_Discrete_Exterior_Calculus_kahler_citations}

\end{document}